\documentstyle{article}

\title{Extrema of $p$-energy functional\\
on a Finsler manifold}
\author{C.Udri\c ste and M.Neagu}
\date{}
\begin{document}
\nocite{3,5,6,7}
\maketitle

\begin{abstract}
In $\S 1$ the authors present the basic properties of the scalar product
along a curve on a Finsler manifold. In
$\S 2$ they investigate the variational formulae for the $p$-energy
functional ($p\in R-\{0\}$). This concept generalises the notions of length
($p=1$) and energy ($p=2$) of a curve.
$\S 3$ analyses the extrema of $p$-energy when
the Finsler space has constant curvature.
\end{abstract}

{\bf \S 1. Preliminaries}
\medskip

Let $(M,F)$ be a connected $n$-dimensional Finsler manifold whose
fundamental function
$F:TM\rightarrow R$ verifies the following axioms:

$\;\;$(F1) $F(x,y)>0;\; \forall x\in M,\; \forall y\ne 0$.

$\;\;$(F2) $F(x,\lambda y)=\vert\lambda\vert F(x,y); \; \forall\lambda\in R,
\forall (x,y)\in TM$.

$\;\;$(F3) the fundamental tensor $\displaystyle{g_{ij}(x,y)={1\over 2}{\partial^2F^2\over\partial
y^i\partial y^j}}$ is positive definite.

$\;\;$(F4) $F$ is $C^\infty$ at every point $(x,y)\in TM$ with $y\ne 0$
and continuous at every $(x,0)\in TM$.
Then, the absolute Finsler energy is $F^2(x,y)=g_{ij}(x,y)y^iy^j$.

Let $c:[a,b] \to M$ be a $C^{\infty}$ regular curve on $M$.
For any two vector fields
$X(t)=X^i(t)\displaystyle{\left.{\partial\over\partial x^i}
\right\vert_{c(t)}},
\; Y(t)=Y^i(t)\displaystyle{\left.{\partial\over\partial x^i}
\right\vert_{c(t)}}$
along the curve $c$,
we introduce \cite {1}, \cite {6}
the {\it scalar product
$
g(X,Y)(c(t))=g_{ij}(c(t),\dot c(t))X^i(t)Y^j(t)
$
along the curve $c$}.

{\bf Remarks.}\\
i) If $X=Y$, then we obtain
$\Vert X\Vert=\sqrt{g(X,X)}.$\\
ii) The vector fields $X$ and $Y$ {\it are orthogonal along the curve $c$}
and we write $X\bot Y$ iff $g(X,Y)=0$.

Let $C\Gamma (N)=(L^i_{jk},N^i_j,C^i_{jk})$ be the {\it Cartan canonical
$N$-linear connection} determined by the fundamental tensor $g_{ij}(x,y)$.
The coefficients of this connection are expressed by
$$
L^i_{jk}=\displaystyle{{1\over 2}g^{im}\left(
{\delta g_{mk}\over\delta x^j}+{\delta g_{jm}\over\delta x^k}-
{\delta g_{jk}\over\delta x^m}\right)},
\;
C^i_{jk}=\displaystyle{{1\over 2}g^{im}\left(
{\partial g_{mk}\over\partial y^j}+{\partial g_{jm}\over\partial y^k}-
{\partial g_{jk}\over\partial y^m}\right)},
$$
$$
N^i_j=\displaystyle{{1\over 2}{\partial\over\partial y^j}\left(
\Gamma^i_{kl}y^ky^l\right)={1\over 2}{\partial\Gamma^i_{00}\over\partial y^j}},
\;
\Gamma^i_{jk}=\displaystyle{{1\over 2}g^{im}\left(
{\partial g_{mk}\over\partial x^j}+{\partial g_{jm}\over\partial x^k}-
{\partial g_{jk}\over\partial x^m}\right)},
$$
where $\displaystyle{{\delta\over\delta x^i}={\partial\over\partial x^i}+
N_i^j{\partial\over\partial y^j}}$.
\medskip

Let $X$ be a vector field along the curve $c$ expressed locally
by $\displaystyle{X(t)=X^i(t)\left.{\partial\over\partial x^i}
\right\vert_{c(t)}}$.
Using the Cartan $N$-linear connection we define {\it
the covariant derivative
along the curve $c$}, by
$$
{\nabla X\over dt}=
\left\{{dX^i\over dt}+X^m\left [L_{mk}^i(c(t),\dot c(t)){dc^k\over dt}+
C_{mk}^i(c(t),\dot c(t)){\delta \over dt}\left({dc^k \over dt}\right)
\right ]\right\}
\left.{\partial\over\partial x^i}\right\vert_{c(t)}.
$$
Since $\displaystyle{{\delta \over dt}\left({dc^k \over dt}\right)
={d^2c^k \over dt^2}+N_l^k(c(t),\dot c(t))
{dc^l\over dt}}$ we obtain
$$
{\nabla X\over dt}=\left\{{dX^i\over dt}+
X^m\left [\Gamma_{mk}^i(c(t),\dot c(t)){dc^k\over dt}+C_{mk}^i(c(t),\dot c(t))
{d^2c^k\over dt^2}\right ]\right\}\left.{\partial\over\partial x^i}
\right\vert_{c(t)},
$$
where $\Gamma_{mk}^i(c(t),\dot c(t))=L_{mk}^i(c(t),\dot c(t))+
C_{ml}^i(c(t),\dot c(t))N_k^l(c(t),\dot c(t))$.

{\bf Remarks.}\\
i) $c$ is a {\it geodesic} iff $\displaystyle{{\nabla\dot c\over dt}=0}$.\\
ii) Since $C\Gamma (N)$ is a metrical connection we have\\
$$
\displaystyle{{d\over dt}\left[g(X,Y)\right]=g\left({\nabla X\over dt},Y\right)+
g\left(X,{\nabla Y\over dt}\right).}
$$

{\bf \S2. Variations and extrema of $p$-energy functional}
\medskip

Let $x_0,x_1\in M$ be two points not necessarily distinct. We denote\\
$\Omega\stackrel{\hbox{not}}{=}
\{c:[0,1]\rightarrow M\vert$ c is piecewise $C^\infty$
regular curve, $c(0)=x_0$, $c(1)=x_1\}$.

For every $p\in R-\{0\}$ we define the $p-energy\; functional$
$$
E_p:\Omega\rightarrow R_+,
$$
$$
E_p(c)=\int_0^1\left [g_{ij}(c(t),\dot c(t))
{dc^i\over dt}{dc^j\over dt}\right]^{p/2}dt=
\int_0^1[g(\dot c,\dot c)]^{p/2}dt=
\int_0^1\Vert\dot c\Vert^pdt.
$$

{\bf Remarks.}\\
i) This general functional was studied for the first time by Udriste
\cite{7}-\cite{10} on Riemannian manifolds. In their papers \cite{2}, \cite{3},
de Cecco and Palmieri study  the same functional
for $p\in(1,\infty)$, but from
a topological point of view, ignoring the geometrical structure.\\
ii) For $p=1$ we obtain the {\it length functional}
$\displaystyle {L(c)=\int_0^1\Vert\dot c\Vert dt}\;$,
and for $p=2$ we obtain the {\it energy functional}
$E(c)=\displaystyle{\int_0^1\Vert\dot c\Vert^2 dt}\;$.\\
iii) For any naturally parametrized curve
(i.e. ,$\; \Vert\dot c\Vert =$ constant)
we have $E_p(c)=$
$=(L(c))^p=(E(c))^{p/2}$.\\
iv) The $p$-energy of a curve is dependent of parametrization if $p\neq 1$.
\medskip

For every $c\in\Omega$ we denote
$T_c\Omega=\{X:[0,1]\rightarrow TM\vert X$ is continuos, piecewise
C$^\infty, X(t)\in T_{c(t)}M,\;\forall t\in [0,1],\;X(0)=X(1)=0\}$.
Let $(c_s)_{s\in (-\epsilon , \epsilon )}\subset\Omega$ be
one parameter variation of the curve $c\in\Omega$.
We denote $X(t)=\displaystyle{{dc_s\over ds}(0,t)\in T_c\Omega}$.
Using the equality
$\displaystyle{g\left({\nabla\dot c_s\over \partial s},\dot c_s\right)=
g\left({\nabla\over\partial t}\left({\partial c_s\over\partial s}\right),\dot c_s\right)}$
we can prove the following

{\bf Theorem.} {\it The first variation of the $p$-energy is

$\displaystyle{{1\over p}{dE_p(c_s)\over ds}(0)=
-\sum\limits_tg(X,\Delta_t(\Vert\dot c
\Vert^{p-2}\dot c))-}$
$$
-\int_0^1\Vert\dot c\Vert^{p-4}
g\left (X,\Vert\dot c\Vert^2
{\nabla\dot c\over dt}+(p-2)g
\left({\nabla\dot c\over dt},\dot c\right)\dot c\right )dt,
$$
where $\Delta_t(\Vert\dot c\Vert^{p-2}\dot c)=
(\Vert\dot c\Vert^{p-2}\dot c)_{t^+}-
(\Vert\dot c\Vert^{p-2}\dot c)_{t^-}$
represents the jump of $\Vert\dot c\Vert^{p-2}\dot c$
at the discontinuity point $t\in (0,1). $}

{\bf Corollary.}
{\it The curve $c$ is a critical point of $E_p$
iff $c$ is a geodesic}.

{\bf Remark.}
For $p=1$ the curve $c$ is a reparametrized geodesic.\medskip

Now, let $c\in\Omega$ be a critical point for $E_p$
(i.e. the curve $c$ is a geodesic).
Let $(c_{s_1s_2})_{s_1,s_2\in(-\epsilon ,\epsilon)}\subset\Omega$ be
a two parameter variation of $c$. Using the notations:\\
$X(t)=\displaystyle{{\partial c_{s_1s_2}\over\partial s_1}}(0,0,t)\in T_c\Omega$,
$Y(t)=\displaystyle{{\partial c_{s_1s_2}\over\partial s_2}}(0,0,t)\in T_c\Omega$,
$\Vert\dot c\Vert =v=$constant and
$I_p(X,Y)=\displaystyle{
{\partial^2E_p(c_{s_1s_2})
\over\partial s_1\partial s_2}}(0,0)$, we obtain the following\medskip

{\bf Theorem.} {\it
The second variation of the $p$-energy is\medskip

$\displaystyle{{1\over pv^{p-4}}}I_p(X,Y)=-
\sum\limits_tg\left (Y,v^2\Delta_t\left(
{\nabla X\over dt}\right)+(p-2)g\left(\Delta_t\left({\nabla X\over dt}
\right),\dot c\right)\dot c\right )-$
$$-\int_0^1g\left(Y,v^2\left [{\nabla\over dt}{\nabla X\over dt}+R^2(X,\dot c)\dot c\right ]+
(p-2)g\left(\left[{\nabla\over dt}{\nabla X\over dt}+
R^2(X,\dot c)\dot c\right],\dot c\right)\dot c\right)dt,$$ where
$\displaystyle{\Delta_t\left ({\nabla X\over dt}\right )=
\left ({\nabla X\over dt}\right )_{t^+}-
\left ({\nabla X\over dt}\right )_{t^-}}$
represents the jump of
$\displaystyle{{\nabla X\over dt}}$
at the discontinuity point $t\in (0,1)$ and, if $R^l_{ijk}(c(t),\dot c(t))$
represents the components of Finsler $h$-curvature, then\\
$$R^2(X,\dot c)\dot c=
R^l_{ijk}(c(t),\dot c(t))\displaystyle{{dc^i\over dt}
{dc^j\over dt}X^k{\partial\over\partial x^l}}=R^l_{jk}(c(t),\dot c(t))
{dc^j\over dt}X^k{\partial\over\partial x^l}.$$
}

{\bf Remark.}
It is well known that we have
$$
R^i_{jk}=\displaystyle{{\delta N^i_j\over\delta x^k}-
{\delta N^i_k\over\delta x^j}},\;
R^i_{hjk}=\displaystyle{{\delta L^i_{hj}\over\delta x^k}-{\delta L^i_{hk}\over\delta x^j}
+L^s_{hj}L^i_{sk}-L^s_{hk}L^i_{sj}+C^i_{hs}R^s_{jk}}.
$$
Moreover, using the Ricci identities for the deflection tensors, we also have
$$
R^i_{jk}=R^i_{mjk}y^m=R^i_{0jk}.
$$
\medskip

{\bf Corollary.}
{\it $I_p(X,Y)=0; \forall\;Y\in T_c\Omega\Leftrightarrow X$
is a Jacobi field in the sense of
Matsumoto (i.e. $\displaystyle{{\nabla\over dt}
{\nabla X\over dt}+R^2(X,\dot c)
\dot c=0}$. See} \cite {6}{\it, pp. $289$).}\\

In these conditions we have the following

{\bf Definition.}
A point $c(b)$, $0\leq a<b<1$, of a geodesic $c\in\Omega$ is called a
{\it conjugate point} of a point $c(a)$ along the curve $c$,
if there exists a non-zero Jacobi field which vanishes at $t\in\{a,\;b\}.$

Now, integrating by parts and using the property of metrical connection we find
$${1\over pv^{p-4}}I_p(X,Y)=\int_0^1\left \{v^2\left [
g\left({\nabla X\over dt},
{\nabla Y\over dt}\right)-R^2(X,\dot c,Y,\dot c)\right ]+\right. $$\\
$$\left. +(p-2)g\left (\dot c,{\nabla X\over dt}\right )g\left (
\dot c,{\nabla Y\over dt}\right )\right \}dt,$$\\
where $R^2(X,\dot c,Y,\dot c)=g(R^2(Y,\dot c)\dot c,X)=
R_{0i0j}(c(t),\dot c(t))X^iY^j$.\medskip

{\bf Remark.}
Let $R_{ijk}=g_{jm}R^m_{ik}$. In any Finsler space
it is satisfied the identity,
$$R_{ijk}+R_{jki}+R_{kij}=0,$$
obtained by the Bianchi identities.
Because $R_{0i0j}=R_{i0j}=R_{j0i}=R_{0j0i}$ we obtain
$R^2(X,\dot c,Y,\dot c)=R^2(Y,\dot c,X,\dot c)$.\medskip

The quadratic form associated to the Hessian of the $p$-energy is

$I_p(X)\stackrel{\hbox {not}}{=}I_p(X,X)=
\displaystyle{\int_0^1\left\{v^2\left [\left\Vert{\nabla X\over dt}
\right\Vert^2-R^2(X,\dot c,X,\dot c)\right ]+\right.}$
$$\left. +(p-2)\left [g\left(\dot c,{\nabla X\over dt}\right)
\right ]^2\right\}dt. $$

{\bf Lemma 1.} {\it
Let $T_c^\bot\Omega=\{X\in T_c\Omega\vert
g(X,\dot c)=0\}$
and $T_c^{'}\Omega=\{X\in T_c\Omega\vert X=f\dot c$,
where $f:[0,1]\rightarrow R$ is continuous, piecewise $C^\infty ,\;
f(0)=f(1)=0\}$. Then\\
i) $T_c\Omega=T_c^\bot\Omega\oplus
T_c^{'}\Omega$;
ii) $I_p(T_c^{\bot}\Omega , T_c^{'}\Omega )=0$.}

{\bf Proof.}\\
ii) Let $X\in T_c^{\bot}\Omega$ and
$Y=f\dot c\in T_c^{'}\Omega. $
Since $X\bot\dot c$ and $c$ is a geodesic, we have
$\displaystyle{g\left(\dot c,{\nabla X\over dt}\right)=0}$.
In these conditions it follows
$${1\over pv^{p-4}}I_p(X,f\dot c)=
\int_0^1\left\{v^2\left [g\left({\nabla X\over dt},f'\dot c\right)-
R^2(X,\dot c,f\dot c,\dot c)\right ]+\right.$$
$$\left.+(p-2)g\left (\dot c,{\nabla X\over dt}\right )
g\left (\dot c,f'\dot c\right )\right\}dt.$$
Hence $I_p(X,f\dot c)=0$. $\Box$
\medskip

{\bf Remark.}
According to the preceding lemma,
the spaces $T^{\bot}_c \Omega $ and $T'_c \Omega $ are orthogonal with respect
to the bilinear form $I_p$ and consequently,
the study of the signature of the quadratic
form $I_p$ is reduced to the study of signatures of its
restrictions to $T'_c\Omega$ and $T_c^{\bot}\Omega. $
\medskip

{\bf Proposition 1.} {\it
Let $c$ be a geodesic and $p\in R-\{0,1\}$. Then\\
i) $I_p(T'_c\Omega )\ge 0$ for $p\in(-\infty ,0)\cup (1,\infty )$,\\
ii) $I_p(T'_c\Omega )\le 0$ for $p\in (0,1)$.\\
Moreover, in both cases: $I_p(X)=0\Leftrightarrow X=0. $}

{\bf Proof.}
Let $X=f\dot c\in T'_c\Omega$. Then we have
\begin{flushleft}
$
\displaystyle{{1\over v^{p-4}}I_p(X)=p\int_0^1\left \{
v^2\left [g(f'\dot c,f'\dot c)-
R^2(f\dot c,\dot c,f\dot c,\dot c)\right ]+
(p-2)\left [g(\dot c,f'\dot c)\right ]^2\right \}dt=}
$\\
$
\displaystyle{=p\int_0^1\left [v^4(f')^2+(p-2)v^4(f')^2\right ]dt=
\int_0^1p(p-1)v^4(f')^2dt}
$.\\
\end{flushleft}
Moreover, if
$I_p(X)=0\Leftrightarrow f'=0
\Leftrightarrow f$ is constant.
The conditions $f(0)=f(1)=0$ imply $f=0.\:\Box$

Because $I_p(T'_c\Omega)$ is positive definite for $p\in(-\infty ,0)\cup
(1,\infty)$ and negative definite for $p\in(0,1)$,
it is sufficient to study the behaviour of $I_p$
restricted to $T_c^{\bot}\Omega$.
Since $X\bot\dot c$ and the curve $c$ is a geodesic it follows
$\displaystyle{g\left(\dot c,{\nabla X\over dt}\right)=0}$.
Hence, for all $X\in T_c^{\bot}\Omega$, we have
$$
{1\over pv^{p-2}}I_p(X)=\int_0^1\left [
\left\Vert {\nabla X\over dt}\right\Vert^2-R^2(X,\dot c,X,\dot c)
\right]dt\stackrel{\hbox{not}}
{=}I(X).
$$

{\bf Lemma 2.}
{\it The following statements are equivalent:\\
i) the curve $c$ has no conjugate points to $x_0=c(0)$,\\
ii) $I\vert_{T_c^{\bot}\Omega}$ is positive definite.}\medskip

{\bf Proof.}
The proof of the lemma follows closely the proof of
Kobayashi for the case of a Riemannian manifold (See \cite {4}, vol 2, pp 72-76).

i)$\Rightarrow $ ii).
Let $J_{c(0)}^{\bot}=\{X:[0,1]\rightarrow TM
\vert \;\;X(t)\in T_{c(t)}M, \;\;\forall t\in [0,1]$, $X$ is Jacobi
field, $X(0)=0, X\bot\dot c\}$. Then
$dim_{\scriptstyle R}J_{c(0)}^{\bot}=n-1$,
where $n=dim\;M$.
Let $\{Y_1, Y_2,..., Y_{n-1}\}$ be a basis in $J_{c(0)}^{\bot}$.
Since the geodesic $c$ has no conjugate point to $x_0=c(0)$ it follows that
$\{Y_1, Y_2,..., Y_{n-1}\}\subset T_c^{\bot}\Omega$ is a
basis for $T_c^{\bot}\Omega.$

Let $X\in T_c^{\bot}\Omega $. There exist the functions
$f_1(t),f_2(t),...,f_{n-1}(t)$ such that $X=$
$=\sum_{i=1}^{n-1}f_iY_i$. We have
$$
\displaystyle{g\left({\nabla X\over dt},{\nabla X\over dt}\right)-
R^2(X,\dot c,X,\dot c)=}
$$
$$
\displaystyle{g\left({\nabla X\over dt},{\nabla X\over dt}\right)-
\sum_if_ig(R^2(Y_i,\dot c)\dot c,X)=g\left({\nabla X\over dt},{\nabla X\over dt}\right)+\sum_if_i
g\left({\nabla \over dt}{\nabla Y_i\over dt},X\right)=}
$$
$$
=\displaystyle{g\left(\sum_if_i'Y_i,\sum_jf'_jY_j\right)+2g\left
(\sum_if'_iY_i,\sum_jf_j{\nabla Y_j\over dt}\right)+}
$$
$$
\displaystyle{+g\left
(\sum_if_i{\nabla Y_i\over dt},
\sum_jf_j{\nabla Y_j\over dt}\right)+
g\left(\sum_if_i{\nabla\over dt}{\nabla Y_i\over dt},\sum_jf_jY_j\right)}.
$$
\newpage
On the other hand, we find
\begin{flushleft}

$\qquad\displaystyle{{d\over dt}\left[g\left(\sum_if_iY_i,\sum_jf_j{\nabla Y_j\over dt}
\right)\right]=g\left(\sum_if_i'Y_i,\sum_jf_j{\nabla Y_j\over dt}\right)+}$

$\qquad\displaystyle{+g\left(\sum_if_i{\nabla Y_i\over dt},\sum_jf_j{\nabla Y_j\over dt}
\right)+g\left(\sum_if_iY_i,\sum_jf_j'{\nabla Y_j\over dt}\right)+}$

$\qquad\displaystyle{+g\left(\sum_if_iY_i,\sum_jf_j{\nabla\over dt}
{\nabla Y_j\over dt}\right).}$
\end{flushleft}
Combining these equalities we obtain
\begin{flushleft}

$\qquad\displaystyle{g\left({\nabla X\over dt},{\nabla X\over dt}\right)-
g(R^2(X,\dot c)\dot c,X)=g\left(\sum_if_i'Y_i,\sum_jf_j'Y_j\right)+}$

$\qquad\displaystyle{+{d\over dt}\left[g\left(\sum_if_iY_i,\sum_jf_j{\nabla Y_j\over
dt}\right)\right]+g\left(\sum_if_i'Y_i,\sum_jf_j{\nabla Y_j\over dt}\right)-}$

$\qquad\displaystyle{-g\left(\sum_if_iY_i,\sum_jf_j'{\nabla Y_j\over dt}\right).}$
\end{flushleft}
Because we have $R^2(X,\dot c,Y,\dot c)=R^2(Y,\dot c,X,\dot c)$
, any two Jacobi fields $X$ and $Y$ such that $X(0)=Y(0)=0$ satisfy
$\displaystyle{g\left(X,{\nabla Y\over dt}\right)=g\left({\nabla X\over dt},
Y\right)}.$ Particularly,
$\displaystyle{g\left(Y_i,{\nabla Y_j\over dt}\right)=
g\left({\nabla Y_i\over dt},Y_j\right).}$ In these conditions we have
\begin{flushleft}
$\qquad\displaystyle{g\left(\sum_if_i'Y_i,\sum_jf_j{\nabla Y_j\over dt}\right)-
g\left(\sum_jf_jY_j,\sum_if_i'{\nabla Y_i\over dt}\right)=0},$
\end{flushleft}
and we obtain
\begin{flushleft}
$\qquad\displaystyle{I(X)=\int_0^1g\left(\sum_if_i'Y_i,\sum_jf_j'Y_j\right)dt+
g\left(\sum_if_iY_i,\sum_jf_j{\nabla Y_j\over dt}\right)_{t=1}=}$

$\qquad\displaystyle{=\int_0^1g\left(\sum_if_i'Y_i,\sum_jf_j'Y_j\right)\geq 0}.$
\end{flushleft}
We have $I(X)=0$ iff $f'_i=0,\;
\forall i=\overline {1,n-1}$ iff $X$ is a Jacobi field.
Since the geodesic $c$ has no conjugate points it follows $X=0$.
\medskip

ii) $\Rightarrow $ i).
We assume that $\exists \;x_{t_0}=c(t_0)$, a point which is conjugate to
$x_0=c(0)$ and $t_0\in (0,1)$. Then $\exists \;Y$ a nonzero Jacobi field
such that $Y(0)=Y(t_0)=0$. Let $U$ be a sufficiently small convex
neighborhood and let $\delta >0$ such that $c(t_0-\delta ),
c(t_0+\delta )\in U$. Then there exists a unique Jacobi field $W$
determined by the boundary values $W(t_0-\delta)=
Y(t_0+\delta )$
and $W(t_0+\delta )=0$. The vector field $X$ is defined along $c$ by
$$
X=\left \{ \begin{array}{lll}
Y\;\quad\;\hbox{from}\;c(0)\; \hbox{to}\;c(t_0-\delta)\\
W\;\quad\hbox{from}\;c(t_0-\delta)\; \hbox{to}\; c(t_0+\delta)\\
0\;\quad\;\;\hbox{from}\;c(t_0+\delta)\;\hbox{to}\;c(1).\\
\end{array} \right.
$$
Denoting
$I_a^b(X)=\displaystyle{\int_a^b\left(\left\Vert {\nabla X
\over dt}\right\Vert^2-R^2
(X,\dot c,X,\dot c)\right)dt}$ and using results from the proof
 i) $\Rightarrow $ ii) we obtain
$0=I_0^{t_0}(Y)=I_0^{t_0-\delta}(Y)+I_{t_0-\delta}^{t_0}(Y)$, since
$Y$ is a Jacobi field. In conclusion:
$I(X)=I(X)-I_0^{t_0}(Y)=I_0^{t_0-\delta}(Y)+
I_{t_0-\delta}^{t_0+\delta}(W)-I_0^{t_0-\delta}(Y)
-I_{t_0-\delta}^{t_0}(Y)=$\\
$=I_{t_0-\delta}^{t_0+\delta}(W)-I_{t_0-\delta}^{t_0}(Y). $
Let
$$
\overline {Y}=\left \{ \begin{array}{ll}
Y\; \hbox{from}\; c(t_0-\delta )\;  \hbox{to}\; c(t_0)\\
0\;\; \hbox{from}\; c(t_0)\; \hbox{to}\;c(t_0+\delta )\\
\end{array} \right.
$$ be a piecewise Jacobi field. Then
$I_{t_0-\delta}^{t_0+\delta}(W)<I_{t_0-\delta}^{t_0+\delta}(\overline {Y})=
I_{t_0-\delta}^{t_0}(Y)$ implies
$I(X)<0, $ contradiction. $\Box $\medskip

With Lemma {\bf 2} and the relation between $I_p$ and
$I$ we have

{\bf Proposition 2.}
{\it Let $c\in\Omega $ be a geodesic and $p\in R-\{0,1\}$.
\begin{flushleft}
i) If $c$ has no conjugate points to $x_0=c(0)$, then\\
$\quad I_p(T_c^{\bot}\Omega )
\geq 0\; \hbox{for}\;p\in (0,1)\cup (1,\infty )$ and
$I_p(T_c^{\bot}\Omega )
\leq 0\; \hbox{for}\;p\in (-\infty ,0)$.\\
$\quad$Moreover, in both cases $I_p(X)=0\Leftrightarrow X=0. $
\end{flushleft}
\begin{flushleft}
ii) If $c$ has conjugate points to $x_0=c(0)$, then\\
$\quad\exists X\in T_c^{\bot}\Omega \;\hbox{such that}\;
I_p(X)<0\; \hbox{for}\; p\in (0,1)\cup (1,\infty )$
and\\
$\quad\exists X\in T_c^{\bot}\Omega \;\hbox{such that}\;
I_p(X)>0\; \hbox{for}\; p\in (-\infty ,0)$.
\end{flushleft}}

Combining the propositions {\bf 1} and {\bf 2} we obtain

{\bf Corollary (extrema of the $p$-energy).}\\
{\it Let $p\in R-\{0,1\}$ and $c\in\Omega $ be a geodesic such that
\begin{flushleft}
A) has no conjugate points to $x_0=c(0)$. Then\\
$\quad$i) $c$ did not even minimize, did not even maximize
$E_p$ for $p\in (-\infty ,0)\cup (0,1);$\\
$\quad$ii) $c$ not maximizes $E_p$ for $p\in (1,\infty );$\\
B) has conjugate points to $x_0=c(0)$. Then\\
$\quad$i) $c$ not maximizes $E_p$ for $p\in (-\infty ,0);$\\
$\quad$ii) $c$ not minimizes $E_p$ for $p\in (0,1);$\\
$\quad$iii) $c$ did not even minimize, did not even maximize $E_p$
for $p\in (1,\infty )$.
\end{flushleft}}

{\bf Remarks.}

1) According to the property (F2) imposed to the Finsler metric, the preceding consequence
is valid replacing $x_0$ with $x_1$ by symmetry.

2) For the case $p\in (1,\infty )$, supposing that exists a $minimal\;
geodesic\gamma\in\Omega $
(i.e. it minimizes the length functional),
then $\gamma $ is a global minimum point
for the $p$-energy $E_p$ since
$E_p(\gamma )=(L(\gamma ))^p
\leq (L(c))^p\leq E_p(c),
\forall c\in\Omega $, where the last inequality is the H\H older
inequality (For details, see \cite {10}).
On the other hand, we have the H\H older inequality for the case $p\in (0,1)$
$$\displaystyle{\int_0^1\vert fg\vert dt\geq
\left (\int_0^1\vert f\vert^pdt\right )^{1\over p}
\left (\int_0^1\vert g\vert^qdt\right )^{1\over q}},$$
$\hbox{where}\;\;
q=\displaystyle{{p\over {p-1}}}\in (-\infty ,0).$
In these conditions it follows
\begin{flushleft}
i)$E_p(c)\leq (L(c))^p\;
\hbox{for}\; p\in (0,1)$\\
ii)$E_p(c)\geq (L(c))^p\;
\hbox{for}\; p\in (-\infty ,0)$, for any curve $c\in\Omega $.
\end{flushleft}
In conclusion we have
\begin{flushleft}
i)$E_p(\gamma )=(L(\gamma )
)^p\leq (L(c))^p\geq E_p(c)$
for $p\in (0,1)$\\
ii)$E_p(\gamma )=(L(\gamma ))^p
\geq (L(c))^p\leq E_p(c)$
for $p\in (-\infty ,0)$.
\end{flushleft}
It follows that, in the cases $p\in (0,1)\cup (-\infty ,0)$, the H\H older formula did
not decide upon the role of minimal geodesics as extremum points of
$E_p$.
Actually, the statement A) of the preceding consequence solves this problem.
\medskip

{\bf \S 3. Extrema of $p$-energy on constant curvature Finsler spaces}
\medskip

We assume the Finsler space ($M$,$F$) is complete, of dimension $n\geq 3$
and of constant curvature $K\in R$. Hence, we have
$$H_{ijkl}=K(g_{ik}g_{jl}-g_{il}g_{jk}),$$
where $H_{ijkl}$ are the components of the $h$-curvature tensor $H$ of the
Berwald connection $ B\Gamma $. It follows that
$$\displaystyle{R_{ijk}=KF\left(g_{ik}{y_j\over F}-g_{ij}{y_k\over F}\right)},$$
where $y_j=g_{jk}y^k$.
We also have
$$R_{i0k}=R_{ijk}y^j=K(g_{ik}F^2-y_iy_k).$$
Hence, along the geodesic $c\in\Omega$, we obtain
$$R^2(X,\dot c)\dot c=K\{\Vert\dot c\Vert^2X-g(X,\dot c)\dot c\}.$$

{\bf Remark.} This equality is also true in the case of constant $h$-curvature
for the Cartan canonical connection.

As in Matsumoto (see \cite{6}, pp 292) we have\\
{\it i) If $K\leq 0$, then the geodesic $c$ has no conjugate points to $x_0=c(0)$.\\
ii) If $K\geq 0$ and the geodesic $c$ has conjugate points to $x_0=c(0)$,
then the number of conjugate points is finite (Morse index theorem
for a Finsler manifold).}

Moreover, in the case ii), choosing an orthonormal frame of vector fields
$\{E_i\}_{i=\overline{1,n-1}}\in
T^{\bot}_c\Omega$
parallelly propagated along the geodesic $c$,
we can build a basis $\{U_i,V_i\}_{i=\overline{1,n-1}}$
in the set of Jacobi fields orthogonal to $\dot c$, defining
$$U_i(t)=\sin(\sqrt Kvt)E_i\;\hbox{and}\;V_i(t)=\cos(\sqrt Kvt)E_i,$$
where $v=\Vert\dot c\Vert=\hbox{constant}$.
In conclusion, the distance between two consecutive
conjugate points is $\displaystyle{{\pi\over{\sqrt K}}}$.In these conditions
we can prove the following

{\bf Theorem.}
{\it Let $(M,F)$ be a Finsler space, as above, and let $c=c_p\in\Omega $
be a global extremum point for the $p$-energy functional $E_p$,
where $p$ is a number in $R-\{0,1\}$.In these conditions
we have\\
i)If $p\in (-\infty ,0)$, then $c$ has conjugate points, $K>0$ and
$$\left [{(m(c)+1)\pi\over\sqrt{K}}\right ]^p\leq
E_p(c)\leq\left [{m(c)\pi\over\sqrt{K}}\right ]^p,\;$$\\
where $m(c)$ is the maximal number of conjugate points to $x_0=c(0)$
along the geodesic $c$.\\
ii) If $p\in (0,1)$, then $c$ has conjugate points, $K>0$ and
$$\left [{m(c)\pi\over\sqrt{K}}\right]^p\leq
E_p(c)\leq\left [{(m(c)+1)\pi\over\sqrt{K}}\right ]^p.\;$$
iii) If $p\in (1,\infty)$, then $c$ is a minimal geodesic
(i.e. it minimizes the length functional)}.

{\bf Proof.}\\
i) If $p\in (-\infty ,0)$ and $c$ is an extremum point for the
$p$-energy $E_p$, then $c$ is a minimum point
and the curve $c$ must have conjugate points
to $x_0$, respectively to $x_1$, and hence $K>0$.
Let $x_0^1,x_0^2,...,x_0^{m(c)}$
be the consecutive conjugate points to $x_0$.
Since the distance between two consecutive conjugate points is
$\displaystyle{{\pi\over\sqrt{K}}}$ it follows
$\displaystyle{{m(c)\pi\over\sqrt{K}}\leq L(c)
\leq {(m(c)+1)\pi\over\sqrt{K}}}$.
On the other hand $E_p(c)=(L(c))^p$, and hence, the above inequality is true.\\
ii) By analogy to i).\\
iii) By the above Remark 2), if $\gamma\in\Omega $ is a minimal geodesic,
then $E_p(\gamma )\leq E_p(c)$.
But $c$ is a minimum point for $E_p$,
and hence $E_p(c)\leq E_p(\gamma )$. In conclusion, we have
$E_p(\gamma )=(L(\gamma ))^p=(L(c))^p=E_p(c)$
and consequently $L(\gamma )=L(c)$.
Hence $c$ is a minimal geodesic.
$\Box $
\medskip

If we denote $m=\sup\{m(c)\;\vert \;c\in\Omega ,\;c-
\hbox{geodesic}\}\in N$, we obtain the following

{\bf Corollary.}{\it If there is $c\in\Omega$ a global extremum point for the
$p$-energy functional $E_p$, where $p\in(-\infty ,0)\cup (0, 1)$,
we must have $m<\infty$ and $m(c)=m$.}
\medskip

{\bf Remarks.}\\
i) If $x_1$ is not a conjugate point to $x_0$, then
it follows
$\displaystyle{\Vert\dot c\Vert ={(m+1)\pi\over \sqrt K}}$ and
$\displaystyle{E_p(c)=\left[{(m+1)\pi\over\sqrt K}\right]^p}$,
because the $p$-energy
is dependent of parametrization.\\
ii) If $x_1$ is a conjugate point to $x_0$, then we obtain
$\displaystyle{E_p(c)=\left({m\pi\over\sqrt K}\right)^p}$ and
$\displaystyle{v={m\pi\over\sqrt K}}$.
\medskip

{\bf One example.} In the case of Riemannian unit sphere
$S^n\subset R^{n+1},\;n\geq 2$,
it is well known that the geodesics are precisely the great circles, that is
the intersections of $S^n$ with the hyperplanes trough the center of $S^n$.
Moreover, two arbitrary points on $S^n$ are conjugate along a geodesic $\gamma$
if they are antipodal points. In these conditions, for any two points $x_0$ and
$x_1$ on the sphere $S^n$, there is no geodesic trough these points
which has a finite maximal number of conjugate points, because
we can surround the sphere infinite times. Hence, for the unit sphere $S^n$,
we have $m=\infty$. In conclusion, in the case $p\in (-\infty ,0)\cup (0,1)$,
the $p$-energy functional on the sphere has no global extremum points.\\
{\bf Aknowledgements.} We have benefitted from the criticisms of Prof.Dr.
Lajos Tamassy and his co-workers upon previous variants of our paper.

\begin{tabbing}
\kill\\
\kill\\
\end{tabbing}
\begin{center}
University POLITEHNICA of Bucharest\\
Department of Mathematics I\\
Splaiul Independentei 313\\
77206 Bucharest, Romania\\
e-mail:udriste@mathem.pub.ro\\
e-mail:mircea@mathem.pub.ro\\
\end{center}

\end{document}